
\input amstex
\input amsppt.sty
\hsize 30pc
\vsize 47pc
\def\nmb#1#2{#2}         
\def\cit#1#2{\ifx#1!\cite{#2}\else#2\fi} 
\def\totoc{}             
\def\idx{}               
\def\ign#1{}             

\define\X{\frak X}
\define\al{\alpha}

\define\ga{\gamma}
\define\de{\delta}

\define\rh{\rho}
\define\si{\sigma}

\define\ph{\varphi}

\define\ps{\psi}
\define\om{\omega}
\define\Ga{\Gamma}

\define\Th{\Theta}
\define\La{\Lambda}

\define\Ph{\Phi}

\define\Om{\Omega}
\redefine\i{^{-1}}
\define\row#1#2#3{#1_{#2},\ldots,#1_{#3}}

\define\sign{\operatorname{sign}}
\redefine\L{\Cal L}
\define\Fl{\operatorname{Fl}}
\def\today{\ifcase\month\or
 January\or February\or March\or April\or May\or June\or
 July\or August\or September\or October\or November\or December\fi
 \space\number\day, \number\year}
\topmatter
\title A common generalization of the Fr\"olicher-Nijenhuis bracket
and the Schouten bracket for symmetric multi vector fields
\endtitle
\author  Michel Dubois-Violette\\
Peter W. Michor  \endauthor
\leftheadtext{\smc M\. Dubois-Violette, P\. Michor}
\rightheadtext{\smc Fr\"olicher-Nijenhuis-Schouten bracket}
\affil
LPTHE Universit\'e Paris XI, B\^atiment 211, F-91405 Orsay Cedex,
France\\
Erwin Schr\"odinger International Institute of Mathematical Physics,
Wien, Austria
\endaffil
\address
M\. Dubois-Violette:
Laboratoire de Physique Th\'eorique et Hautes Energies,
Universit\'e Paris XI, B\^atiment 211, F-91405 Orsay Cedex,
France
\endaddress
\email flad\@qcd.th.u-psud.fr \endemail
\address
P\. Michor: Institut f\"ur Mathematik, Universit\"at Wien,
Strudlhofgasse 4, A-1090 Wien, Austria; and:
Erwin Schr\"odinger International Institute of Mathematical Physics,
Pasteurgasse 6/7, A-1090 Wien, Austria
\endaddress
\email Peter.Michor\@esi.ac.at \endemail
\date {\today} \enddate
\keywords Fr\"olicher-Nijenhuis bracket, symmetric Schouten bracket,
Poisson bracket \endkeywords
\subjclass 58A10, \endsubjclass
\abstract There is a canonical mapping from the space of sections of
the bundle $\La T^*M\otimes STM$ to $\Om(T^*M;T(T^*M))$. It is shown
that this is a homomorphism on $\Om(M;TM)$ for the
Fr\"olicher-Nijenhuis brackets, and also on $\Ga(STM)$ for the
Schouten bracket of symmetric multi vector fields. But the whole
image is not a subalgebra for the Fr\"olicher-Nijenhuis bracket on
$\Om(T^*M;T(T^*M))$.
\endabstract
\endtopmatter

\document

\heading Table of contents \endheading
\noindent 1. Introduction \leaders \hbox to 1em{\hss .\hss }\hfill {\eightrm
1}\par
\noindent 2. The Poisson bracket for differential forms \leaders \hbox to
1em{\hss .\hss }\hfill {\eightrm 3}\par
\noindent 3. The Fr{\accent "7F o}licher-Nijenhuis bracket on $\Omega
(T^*M;T(T^*M))$ \leaders \hbox to 1em{\hss .\hss }\hfill {\eightrm 5}\par

\head\totoc\nmb0{1}. Introduction \endhead

It is well known that there are several extensions of the
bracket of vector fields on a smooth manifold $M$. In
particular, the Fr\"olicher-Nijenhuis bracket extends
the bracket of vector fields to all vector valued
differential forms on $M$, i.e. to $\Om(M;TM)$.
Another classical extension is the Schouten bracket, this
is an extension of the bracket of vector fields to all
symmetric multivector fields, i.e. to $\Ga(STM)$.
The Schouten bracket has a natural interpretation in
terms of Poisson bracket. Indeed, there is an obvious
isomorphism $\pi^*$ of the algebra $\Ga(STM)$ on the
algebra of smooth functions on $T^* M$ which are
polynomial on the fiber. On the other hand there is a
natural symplectic structure on $T^* M$ and the
Schouten bracket corresponds just to the Poisson bracket
under the above isomorphism.

It is very natural, and it is the aim of this paper, to
try to find a common generalization of the two above
brackets.

Let us give an example of problem where such an extension
could be welcome. Suppose that $M$ is equipped with a
riemannian metric $g$ and let $\underline g$ denote the corresponding
contravariant symmetric two-tensor field. Then,
$\pi^*(\underline g)$ is the hamiltonian of the
geodesic flow on $T^* M$ and the symmetric tensor
fields $S$ satisfying $[\underline g,S]=0$ correspond to functions
on $T^* M$ which are invariant by the geodesic flow;
such symmetric tensor fields are called Killing tensors.
These Killing tensors form a Poisson subalgebra of
$\Ga(STM)$. Now if $S$ is a Killing tensor of order
$k$; then it is not hard to show that only its covariant
derivatives of order not greater than $k$ are independent,
i.e. its covariant derivatives of order greater than $k$
are linear combination of the one of order not greater
than $k$ with coefficients which are covariant expressions
in the curvature tensor. This implies in particular that
the equations $[\underline g,S]=0$ have a lot of integrability
conditions and, since these integrability conditions are
always consequence of $d^2=0$, it is natural to introduce
the algebra $\Om(M;STM)$ of symmetric multivector
valued forms to analyse them. This algebra is a
graded-commutative algebra for the graduation given by the
form-degree and on this algebra there is a nice algebra of
graded derivations associated with the metric. Is is
generated by three antiderivations, $\nabla,
\de_g,\de'_g$, where $\nabla$ is the exterior
covariant differential corresponding to the Levi-Civita
connection, $\de_g$ is the unique $C^\infty(M)$-linear
antiderivation such that $\de_g X \in \Om^1(M)$ for
$X\in \Ga(TM)$ is the one-form
$Y\mapsto\de_gX(Y)=g(X,Y)$ and
$\de_g\Om^1(M)=0$, $\de'_g$ is the unique
$C^\infty(M)$-linear antiderivation such that $\de'_g
\Ga(TM)=0$ and $\de'_g\om\in \Ga(TM)$ for
$\om\in\Om^1(M)$ is the vector field obtained by the
contraction of $\underline g$ with $\om$. One has: $\de^2_g=0$,
${\de'}^2_g=0$, $\de_g\de'_g+\de'_g\de_g$ equals the
total degree in form and tensor,
$\nabla\de_g+\de_g\nabla=0$ (because $\nabla$ is
torsion free) and the derivation
$D=\nabla\de'_g+\de'_g\nabla$ is an extension to
$\Om(M;STM)$ of the Schouten bracket with $\underline g$. So it
is natural to try to construct a bracket on
$\Om(M;STM)$ extending the Schouten bracket for which
$D$ is the bracket with $\underline g$. It is not
difficult to construct such a bracket namely $$[\alpha
\otimes F, \beta\otimes G]_\nabla =
L^\nabla_{\alpha\otimes F}(\beta) G - (-1)^{ab}
L^\nabla_{\beta\otimes G} (\alpha) F + \alpha\wedge \beta
\otimes [F,G]$$ for $\alpha\in \Om^a(M), \
\beta\in\Om^b(M),\ F,\ G\in \Ga(STM)$ with
$$
L^\nabla_{\alpha\otimes F}(\om)=i_{\alpha\otimes F}
\nabla\om + (-1)^a\nabla i_{\alpha\otimes F}\om
$$
for $\om\in\Om(M)$ and where the generalised
insertion $i$ is defined by
$$
i_{\alpha\otimes
X_{1}\vee\cdots \vee X_k}(\om)=\sum^k_{r=1} \alpha
\wedge i_{X_r}(\om) \otimes X_1\vee
\dots\widehat{X_r}\dots\vee X_r
$$
(the hat meaning omission of this element).

More generally if $\nabla$ is any torsion free linear
connection on $M$, the above formula defines a bracket
$[\quad,\quad]_\nabla$ which is an extension to
$\Om(M;STM)$ of both the Schouten bracket and the
Fr\"olicher-Nijenhuis bracket. Furthermore, this bracket
is a graded derivation in each variable, it is also
graded antisymmetric but unfortunately it does not
satisfy the graded Jacobi identity.

In this paper we shall follow another way: we first send
$\Om(M;STM)$ in $\Om(T^* M)$ by using the
isomorphism $\pi^*$, then we use a construction introduced by
one of us \cit!{5} to send it in $\Om(T^* M;T(T^* M))$
in which there is the Fr\"olicher-Nijenhuis bracket and
we show that this gives injective homomorphisms of graded
Lie algebras for the Fr\"olicher-Nijenhuis bracket on
$\Om(M;TM)$ and the Schouten bracket on
$\Ga(STM)$. But the common generalization of these two brackets
does not exist on the space $\Om(M;STM)$, only on
$\Om(T^*M;T(T^*M))$. This is similar to the common generalization
found by Vinogradov \cit!{14,1} of the Fr\"olicher-Nijenhuis bracket
and the skew symmetric Schouten bracket on $\Ga(\La TM)$, which
exist only on a quotient of a certain space of `superdifferential
operators' on $\Om(M)$.

\head\totoc\nmb0{2}. The Poisson bracket for differential forms
\endhead

\subhead\nmb.{2.1}. Fr\"olicher-Nijenhuis bracket \endsubhead
For the convenience of the reader we review here the theory of graded
derivations of the graded commutative algebra of differential form on
a smooth manifold $M$. See \cit!{2} and \cit!{3} for the original
source, and \cit!{6} or \cit!{4}, sections~8--11, as a convenient
reference, whose notation we follow here.

The space
$\operatorname{Der}(\Om(M))$ of all graded derivations of the graded
commutative algebra of differential forms on $M$ is a graded Lie
algebra with the graded commutator as bracket.
In the following formulas we will always assume that
$K \in \Om^{k}(M;TM)=\Ga(\La^kT^*M\otimes TM)$,
$L\in \Om^l(M;TM)$, $\om\in\Om^q(M)$.
The formula
$$
(i_K\om)(\row X1{k+q-1})
= {\tsize\frac1{k!(q-1)!}} \sum_{\si}
\sign(\si) \om(K(X_{\si1},\dots,X_{\si k}), X_{\si(k+1)},\dotsc)
$$
for $X_i\in \X(M)$ (or $T_xM$) defines an
graded derivation $i_K \in \operatorname{Der}_{k-1}\Om(M)$ and any
derivation $D$ with $D|\Om^0(M)=0$ is of this form.
On $\Om^{*+1}(M,TM)$ (with the grading $*$)
a graded Lie bracket is given by
$[K,L]^{\wedge} = i_KL - (-1)^{(k-1)(l-1)}i_LK $
where $i_K(\ps\otimes X):=i_K(\ps)\otimes X$,
which satifies $i([K,L]^{\wedge}):= [i_K,i_L]$.
It is called the \idx{\it Nijenhuis-Richardson bracket}, see \cit!{11}.

The exterior derivative $d$ is an element of
$\operatorname{Der}_1\Om(M)$. We define
the \idx{\it Lie derivation} $\L_K = \L(K) \in
\operatorname{Der}_k\Om(M)$ by $\L_K := [i_K,d]$.
For any graded derivation $D \in \operatorname{Der}_k\Om(M)$
there are unique $K \in \Om^k(M;TM)$ and $L \in \Om^{k+1}(M;TM)$
such that $D = \L_K + i_L.$
We have $L=0$ if and only if $[D,d]=0$, and $D|\Om^0(M)=0$ if and
only if $K=0$. Clearly $[[\L_K,\L_L],d] =0$, so we have
$[\L(K),\L(L)] = \L([K,L])$
for a uniquely defined $[K,L] \in \Om^{k+l}(M;TM)$. This
vector valued form $[K,L]$ is called the
\idx{\it Fr\"olicher-Nijenhuis bracket} of $K$ and $L$.
It is well behaved with respect to the obvious relation of
$f$-relatedness of tangent bundle valued differential forms, where
$f:M\to N$. For $k=l=0$
it coincides with the Lie bracket of vector fields.
Let the degree of $\om$ be $q$, of $\ph$ be
$k$, and of $\ps$ be $\ell$. Let the other degrees be as
indicated. Then the following formulas hold
$$\allowdisplaybreaks\align
[\L_K,i_L] &= i([K,L]) - (-1)^{k\ell}\L(i_LK)\tag1\\
i(\om\wedge L) &= \om\wedge i(L). \tag2 \\
\L(\om \wedge K)&= \om \wedge \L_K - (-1)^{q+k-1}i(d\om\wedge K). \tag3 \\
[\om\wedge K_1,K_2] &= \om \wedge [K_1,K_2] - (-1)^{(q+k_1)k_2}
     \L(K_2)\om \wedge K_1 \tag4 \\
&\quad +(-1)^{q+k_1}d\om \wedge i(K_1)K_2. \\
[\ph\otimes X, \ps\otimes Y] &= \ph \wedge \ps \otimes [X,Y]
     + \ph \wedge \L_X\ps\otimes Y - \L_Y\ph \wedge \ps \otimes X \tag5\\
&\qquad + (-1)^k\left(d\ph \wedge i_X\ps \otimes Y + i_Y\ph
     \wedge  d\ps \otimes X \right).
\endalign $$

\subhead\nmb.{2.2}. Poisson manifolds
\endsubhead
Let $(M,\rh)$ be a Poisson manifold, that is a smooth manifold $M$
together with a 2-field $\rh\in \Ga(\La^2 TM)$ satisfying
$[\rh,\rh]=0$, where $[\quad,\quad]$ is the Schouten-Nijenhuis
bracket on $\Ga(\La^{*-1}TM)$, see \cit!{7} and \cit!{12}.
Then $\rh$ induces a skew symmetric differential concomitant on
$C^\infty(M,\Bbb R)$ given by $\{f,g\}_\rh = \rh(df,dg)$. The Jacobi
identity for this bracket is equivalent to $[\rh,\rh]=0$, see
\cit!{7}, 1.4 for a nice proof. Here we view $\rh$ as a skew
symmetric biliear form on $T^*M$, but also as a vector bundle
homomorphism $\rh:T^*M\to TM$.

It is well known that for a symplectic manifold $(M,\om)$ with
associated Poisson structure $\rh=\om\i:T^*M\to TM$ we have the
following exact sequence of Lie algebras:
$$0\to  H^0(M)\to C^\infty(M,\Bbb R)  @>H>> \X_\om(M) @>\ga>>
     H^1(M)\to 0\tag1$$
Here $H^*(M)$ is the real De Rham cohomology of $M$, the space
$C^\infty(M,\Bbb R)$ is equipped with the Poisson bracket
$\{\quad,\quad\}_\rh$,
$\X_\om(M)$ consists of all vector fields $\xi$ with $\L_\xi\om=0$
(the locally Hamiltonian vector fields), which is a Lie algebra
for the Lie bracket. Also $H_f$ is the Hamiltonian vector field
for $f\in C^\infty(M,\Bbb R)$ given by $H_f=\rh(df)$, and $\ga(\xi)$
is the cohomology class of $i_\xi\om$. The spaces $H^0(M)$ and
$H^1(M)$ are equipped with the zero bracket.

\subhead\nmb.{2.3}. The graded Poisson bracket for differential forms
\endsubhead
In \cit!{5} the exact sequence \nmb!{2.2}.\thetag1 has been
generalized in the following way. It was stated there for symplectic
manifolds, but the proofs there work without any change also for
Poisson manifolds.

We consider first the space $\Om(M;TM)=
\Ga(\La^*T^*M\otimes TM)$ of tangent bundle valued
differential forms on $M$, equipped with the Fr\"olicher-Nijenhuis
bracket $[\quad,\quad]$.
We first extend $\rh:T^*M\to TM$ to a
module valued graded derivation of degree $-1$ by
$$\align
\rh:\Om(M)&\to \Om(M;TM),\tag1\\
\rh|\Om^0(M)&=0,\quad\text{ and for }\ph_i\in\Om^1(M)\text{ by}\\
\rh(\ph_1\wedge \dots\wedge \ph_k) &= \sum_{i=1}^k(-1)^{i-1}
     \ph_1\wedge \dots\widehat{\ph_i}
     \dots\wedge \ph_k\otimes \rh(\ph_i).
\endalign$$
Then we have the \idx{\it Hamiltonian mapping}
$$\align
H:\Om(M)&\to \Om(M;TM),\tag2\\
H(\ps):&= \rh(d\ps),\\
H(f_0\,df_1\wedge \dots\wedge df_k) &= \sum_{i=0}^k(-1)^{i}
     df_0\wedge \dots\widehat{df_i}
     \dots\wedge df_k\otimes H_{f_i}.
\endalign$$

\proclaim{Theorem}\cit!{5}. Let $(M,\rh)$ be a Poisson manifold. Then on the
space $\Om(M)/B(M)$ of differential forms modulo exact forms there
exists a unique graded Lie bracket $\{\quad,\quad\}^1_\rh$, which is given by
the quotient modulo $B(M)$ of
$$\align
&\{\ph,\ps\}^1_\rh = i(H_\ph)d\ps, \qquad\text{ or}\\
&\{f_0\,df_1\wedge \dots\wedge df_k,g_0\,
     dg_1\wedge \dots\wedge dg_l\}^1_\rh =\tag3\\
&\quad= \sum_{i,j}(-1)^{i+j}\{f_i,g_j\}_\rh\,
     df_0\wedge \dots\widehat{df_i}\dots\wedge df_k\wedge
     dg_0\wedge \dots\widehat{dg_j}\dots\wedge dg_k,
\endalign$$
such that $H:\Om(M)/B(M)\to \Om(M;TM)$ is a homomorphism of graded
Lie algebras.

If $\rh=\om\i$ for a symplectic structure $\om$ on $M$ then we have a
short exact sequence of vector spaces
$$\align
0\to  H^*(M)\to \Om(M)/B(M)  @>H>> &\Om_{\L\om=0}(M;TM) \to \tag4\\
&\to  H^{*+1}(M)\oplus \Ga(E_\om)\to 0
\endalign$$
where $\Ga(E_\om)$ is a space of sections of a certain vector
bundle and where the space $\Om_{\L\om=0}(M;TM)$ is the graded Lie
subalgebra of all $K\in\Om(M;TM)$ such that for the Lie derivative
along $K$ we have $\L_K\om=0$. We also have the exact sequence of
graded Lie algebras
$$0\to  H^*(M)\to \Om(M)/B(M)  @>H>> \Om_\om(M;TM) @>>>
     H^{*+1}(M)\to 0\tag5$$
where now $\Om_\om(M;TM)$ is the graded Lie
subalgebra of all $K\in\Om^k(M;TM)$ such that for the Lie derivative
along $K$ we have $\L_K\om=0$ and
$K+\frac{(-1)^{k+1}}{k+1}\rh(i_K\om)=0$,
and where on the De~Rham cohomology spaces we put the brackets 0.
\endproclaim

See \cit!{5} for the proof of this theorem and for more information.
The step from the sequence \thetag 4 to \thetag5 was noticed in
\cit!{8}. Parts of this theorem were reproved by a different method
in \cit!{1}. We just note here that on $\Om(M)$ itself the bracket
$\{\quad,\quad\}^1_\rh$ is graded anticommutative, but does not satisfy
the graded Jacobi identity, whereas a second form, $\{\ph,\ps\}^2_\rh
=\L_{H(\ph)}\ps$, satisfies the graded Jacobi identity but is not
graded anticommutative, and they differ by something exact.

\head\totoc\nmb0{3}. The Fr\"olicher-Nijenhuis bracket on $\Om(T^*M;T(T^*M))$
\endhead

\subhead\nmb.{3.1} \endsubhead
Let $M$ be a smooth manifold. We consider the cotangent bundle
$\pi:T^*M\to M$, the Liouville form $\Th_M\in \Om^1(T^*M)$, given by
$\Th_M(\xi)=\langle \pi_{T^*M}\xi,T(\pi_M).\xi\rangle_{TM}$, and the
canonical symplectic form $\om_M =-d\Th_M$.

The space $\Ga(STM)$ of symmetric contravariant tensorfields
carries a natural differential concomitant which was found by
Schouten \cit!{13} and which for $X_i,Y_j \in\X(M)$ and for
$f,g\in C^\infty(M,\Bbb R)$ is given by (see \cit!{7})
$$\align
[f,g] &= 0\tag1\\
[X_1\vee&\dots\vee X_k,Y_1\vee\dots\vee Y_l] = \\
&= \sum_{i,j}[X_i,Y_j]\vee X_1\vee\dots\widehat{X_i}\dots\vee X_k\vee
     Y_1\vee\dots\widehat{Y_j}\dots\vee Y_l,\\
[f,Y_1&\vee\dots\vee Y_l]
     = \sum_jdf(Y_j).Y_1\vee\dots\widehat{Y_j}\dots\vee Y_l,\\
\endalign$$
Obviously $\Ga(S^{*+1}TM)$ is a Lie algebra (with grading $*$,
but not a graded Lie algebra).
Any symmetric multivector field $U\in \Ga(S^kTM)$ may be viewed
as a function on $T^*M$ which is homogeneous of degree $k$ on each
fiber. So we have a linear injective mapping
$$\align
\pi^*: \Ga(S^kTM) &\to C^\infty(T^*M,\Bbb R)\\
(\pi^*U)(\ph) &= \langle \ph^k,U\rangle_{TM}
\endalign$$
It is well known that $\pi^*$ is a homomorphism of Lie algebras,
where on $C^\infty(T^*M,\Bbb R)$ we consider the canonical Poisson
bracket $\{\quad,\quad\}$ induced by $\rh=\om_M\i$. See also
\nmb!{3.5}.\therosteritem2.

\subhead\nmb.{3.2} \endsubhead
We consider the pullback $\pi^*:\Om(M)\to \Om(T^*M)$, and we extend
it to the linear mapping
$$\align
&\pi^*:\Ga(\La^kT^*M\otimes S^lTM) \to \Om^k(T^*M),\\
&(\pi^*A)_\ph(\xi_1,\dots,\xi_k) =
     \langle \ph\vee\dots\vee\ph,
     A(T\pi.\xi_1,\dots,T\pi.\xi_k)\rangle_{TM}.\tag1
\endalign$$
The space
$\Ga(\La T^*M\otimes STM)
     =\bigoplus_{k,l}\Ga(\La^kT^*M\otimes S^lTM)$
is a graded commutative algebra with respect to the degree $k$, and
$\pi^*:\Ga(\La T^*M\otimes STM)\to \Om(T^*M)$ is obviously
a homomorphism with respect to the `wedge' products.
In the following we will have always write $\pi^*$ in
front of any tensor field on $M$ which contains vector field
components, but we will suppress it if we consider pullbacks of
functions or differential forms to $T^*M$.

\proclaim{Lemma}
\roster
\item [2] For each $k\ge 0$ and for $l>0$ the mapping
$$h:\Ga(\La^kT^*M\otimes S^lTM)
     @>{\pi^*}>> \Om^k(T^*M) \to \frac{\Om^k(T^*M)}{B^k(T^*M)}
     @>{H}>> \Om^k(T^*M;T(T^*M))$$
     is injective.
\item For $l=0$ the mapping
$$\Om(M) @>{\pi^*}>> \Om(T^*M) \to \frac{\Om(T^*M)}{B(T^*M)}$$
     induces an injective linear mapping
$$\frac{\Om(M)}{B(M)}\to \frac{\Om(T^*M)}{B(T^*M)}.$$
\item Let $I\in\X(T^*M)$ be the vertical homothetic vector field on
     $T^*M$, given by $I(\ph)=\frac{\partial}{\partial t}|_1t\ph$.
     Then for each $l\ge0$ the image of the linear mapping
$$\pi^*:\Ga(\La^kT^*M\otimes S^lTM) \to \Om^k(T^*M)$$
     is the subspace consisting of all horizontal differential forms
     $\Ph\in\Om(T^*M)$ which satisfy $\L_I\Ph=l.\Ph$
\endroster
\endproclaim

\demo{Proof}
Since $\pi:T^*M\to M$ is a homotopy equivalence with homotopy inverse
the zero section, the pullback operator induces an injective linear
mapping $\pi^*:\Om(M)/B(M)\to \Om(T^*M)/B(T^*M)$. This proves
\therosteritem3.

Now let $0\ne A\in \Ga(\La^kT^*M\otimes S^lTM)$. We
consider the vertical vector field $I\in\X(T^*M)$,
$I(\ph)=vl(\ph,\ph)=\frac{\partial}{\partial t}|_1t\ph$. The flow
of $I$ is given by the vertical homotheties $\Fl^I_t(\ph)=e^t\ph$, we
have $(\Fl^I_t)^*\pi^*A = e^{lt}\pi^*A$, and thus
$$i_Id\pi^*A + 0 = \L_I\pi^*A = ddt|_0(\Fl^I_t)^*\pi^*A
     = ddt|_0e^{lt}\pi^*A = l\pi^*A$$
which is not 0 for $l>0$.
Since $\rh:\Om^{>0}(T^*M)\to\Om(T^*M;T(T^*M))$ is injective,
\therosteritem2 follows.

We also conclude the inclusion $\subseteq$ in \therosteritem4. Since
the assertion is local on $M$, for the converse inclusion $\supseteq$
we may use local coordinates on $T^*Q\subset$ as in the beginning of
the proof of lemma \nmb!{3.3}. Then
$I|Q=\sum p_i\frac{\partial}{\partial p_i}$ and any horizontal form
is a sum of expressions like
$\Ph=f(q,p)dq^{i_1}\wedge \dots\wedge dq^{i_p}\in\Om^p(T^*Q)$. Then
$\L_I\Ph=l.\Ph$ means $\L_If=l.f$ from which we conclude that in
multi-index notation we have
$f(q,p)=\sum_{|\al|=l}f_\al(q)p^\al$, which implies the result, since
we use a partition of unity on $M$.
\qed\enddemo

\proclaim{\nmb.{3.3}. Lemma. Collection of formulas} In the following
$X,Y\in\X(M)$ are vector fields, $\ph\in\Om^p(M)$, $\ps\in\Om^q(M)$,
$K\in\Om^k(M;TM)$, $L\in\Om^l(M;TM)$, and $f\in \Om^0(M)$.  Then the
following formulas hold on $T^*M$. We drop $\pi^*$ in front of
pullbacks of differential forms.
\roster
\item $[hX,hY]=h[X,Y]$.
\item $[\rh\ph,\rh\ps]=0$, thus also
     $[h\ph,\rh\ps]=[\rh d\ps,\rh\ph]=0$, etc.
\item $[hX,\rh\ph]=\rh\L_X\ph$, so also
     $[hX,h\ph]=[hX,\rh d\ph]=h\L_X\ph$.
\item $i_{\rh\ph}\ps=0$ and $i_{\rh\ph}\rh\ps =0$, so also
     $\L_{\rh\ph}\ps = 0$, etc.
\item $\L_{\rh\ph}\pi^*X = - i_X\ph$, so also
     $\L_{h\ph}\pi^*X = -i_Xd\ph$.
\item $\L_{hK}f=\L_Kf$, so also $\L_{hK}\ph=\L_K\ph$. Similarly
       $i_{hK}\ph=i_K\ph$.
\item $[hL,hf]=h\L_Lf$.
\item $\L_{hK}\pi^*L = \pi^*[K,L] + (-1)^{(k-1)l}d\pi^*(i_LK)$.
\item $d\L_{hK}\pi^*L =(-1)^k\L_{hK}d\pi^*L= d\pi^*[K,L]$.
\item $i_{\rh\pi^*K}\ps=0$, so also $\L_{\rh\pi^*K}\ps=0$.
\item $\L_{\rh\pi^*K}\pi^*L = -(-1)^{(k-1)l}\pi^*i_LK$.
\item $i_{hK}\pi^*L = \pi^*i_KL$.
\item $i_{hK}d\pi^*L= \pi^*[K,L]
       -(-1)^kd\pi^*(i_KL+(-1)^{(k-1)(l-1)}i_LK)$.
\item $i_{h(X}\rh\ps = -\rh i_X\ps$.
\item $\L_{h\ph}\pi^*L = -(-1)^{pl}i_Ld\ph$.
\item $[\rh\pi^*K,h\ps]=\rh(i_Kd\ps)-(-1)^ki_{hK}h\ps$.
\endroster
\endproclaim

\demo{Proof}
Let us fix local coordinates $q^1,\dots,q^m$ on an open subset $Q$ of
$M$ and induced coordinates $q^i,p_j$ on $T^*Q\subset T^*M$, so that
the Liouville form $\Th|T^*Q=\sum p_idq^i$ and the symplectic form is
given by $\om=-d\Th=\sum dq^i\wedge dp_i$.
We have
$$\alignat2
\om(\tfrac{\partial}{\partial q^i}) &= dp_i &\qquad
     \rh(dp_i)&= \tfrac{\partial}{\partial q^i}\\
\om(\tfrac{\partial}{\partial p_i}) &= -dq^i &\qquad
     \rh(dq^i)&= -\tfrac{\partial}{\partial p_i}
\endalignat$$
so that for $f\in C^\infty(M,\Bbb R)$, $\ph\in \Om^p(M)$, and
$X\in \X(M)$ we get the following local formulas on
$T^*Q\subset T^*M$:
$$\align
hf &= \rh(df) = -\sum \tfrac{\partial f}{\partial q^i}
     \tfrac{\partial}{\partial p_i} \\
h\ph &= \rh(\sum \tfrac{\partial \ph_{i_1\dots i_p}}{\partial q^i}
     dq^i \wedge dq^{i_1}\wedge \dots \wedge dq^{i_p})\\
&= \sum \tilde\ph_{j_1\dots j_p,j}(q) dq^{j_1}\wedge \dots \wedge dq^{j_p}
     \otimes \tfrac{\partial}{\partial p_j} \\
hX &= -\sum \tfrac{\partial X^i}{\partial q^m}
     p_i\tfrac{\partial}{\partial p_m}
     + \sum X^i \tfrac{\partial}{\partial q^i}.
\endalign$$
 From this \therosteritem1 and \therosteritem2 follow by
straightforward computation, whereas \therosteritem3 follows from
contemplating \nmb!{2.1}.\thetag1.

\therosteritem3 then can be proved as follows:
$$\align
[hX,\rh&(f_0df_1\wedge \dots \wedge df_p)] = \L_{hX}\left(
     \sum_i(-1)^{i-1}f_0df_1\wedge \dots \wedge \rh df_i
     \wedge \dots \wedge df_p\right)\\
&= \L_{hX}f_0.\sum_i(-1)^{i-1}df_1\wedge \dots \wedge
     hf_i\wedge \dots \wedge df_p\\
&\quad+ \sum_{1\le j<i}(-1)^{i-1} f_0.df_1\wedge \dots \wedge
     \L_{hX}df_j\wedge \dots\wedge  hf_i\wedge \dots \wedge df_p\\
&\quad+ \sum_i(-1)^{i-1}f_0.df_1\wedge \dots \wedge
     \L_{hX}hf_i\wedge \dots \wedge df_p\\
&\quad+ \sum_{1\le i<j}(-1)^{i-1} f_0.df_1\wedge \dots \wedge
     hf_i\wedge \dots \wedge \L_{hX}df_j\wedge \dots \wedge df_p\\
&= \rh(\L_X (f_0df_1\wedge \dots df_p)),
\endalign$$
where we also use the following special cases of \therosteritem3, which
are immediate from the local formulas:
$$\align
\L_{hX}f &= \left(-\sum \tfrac{\partial X^k}{\partial q^m}
     p_k\tfrac{\partial}{\partial p_m}
     + \sum X^i \tfrac{\partial}{\partial q^i}\right)f =\L_Xf\\
\L_{hX}df &= d\L_{hX}f = d\L_Xf =\L_Xdf\\
\L_{hX}hf &= \left[-\sum \tfrac{\partial X^k}{\partial q^m}
     p_k\tfrac{\partial}{\partial p_m}
     + \sum X^i \tfrac{\partial}{\partial q^i},
     -\sum \tfrac{\partial f}{\partial q^i}
     \tfrac{\partial}{\partial p_i}\right] \\
&= h\L_Xf = \rh\L_Xdf.
\endalign$$

\therosteritem5 is seen as follows:
$$\align
\L(\rh(f_0&df_1\wedge \dots \wedge df_p))\pi^*X
     = i(\rh(f_0df_1\wedge \dots \wedge df_p))d\pi^*X +0\\
&= i\left(\sum_i(-1)^{i-1}f_0df_1\wedge \dots \widehat{df_i}
     \dots \wedge df_p \otimes hf_i\right)d\pi^*X \\
&= -\sum_i(-1)^{i-1}f_0df_1\wedge \dots \wedge i_Xdf_i \wedge
     \dots \wedge df_p
     = -i_X(f_0df_1\wedge \dots \wedge df_p)
\endalign$$
where we use the special case
$$\align
i(hf)d\pi^*X &= i\left(-\sum \tfrac{\partial f_i}{\partial q^j}
     \tfrac{\partial}{\partial p_j}\right)
     \left(\sum \tfrac{\partial X^k}{\partial q^m}p_kdq^m
     + \sum X^kdp_k\right)\\
&= -\L_Xf = -i_Xdf.
\endalign$$

For the proof of the remaining formulas we assume that
$K=\ph\otimes X$ for $\ph\in\Om^k(M)$ with $d\ph=0$, and
$L=\ps\otimes Y$ for $\ps\in\Om^l(M)$ with $d\ps=0$, where
$X,Y\in\X(M)$. We may do this since locally $\Om(M;TM)$ is linearly
generated by such elements. We will use the formulas of \nmb!{2.1} without
explicitly mentioning them. Under this assumptions we have
$$\align
h(\ph\otimes X)&= \rh d\pi^*(\ph\otimes X)
     = -d\pi^*X\wedge \rh\ph + \ph\wedge hX\\
\L(\ph\otimes X) &= \ph\wedge \L_X.
\endalign$$

\therosteritem6 follows from \therosteritem4 via
$$\align
\L_{hK}f &= i_{hK}df = i\left(-d\pi^*X\wedge\rh\ph
     + \ph \wedge hX \right)df\\
&= -d\pi^*X\wedge i_{\rh\ph}df + \ph \wedge i_{hX}df \\
&= 0 + \ph \wedge i_Xdf = i_Kdf = \L_Kf
\endalign$$
Then we get in turn
$$\align
i_{hK}(f_0df_1\wedge \dots \wedge df_p)
     &= \sum_i(-1)^{i-1}f_0df_1\wedge \dots i_{hK}df_i
     \dots \wedge df_p \otimes hf_i \\
&= i_K(f_0df_1\wedge \dots \wedge df_p),\\
\L_{hK}df&=(-1)^kd\L_{hK}f =(-1)^kd\L_Kf =\L_Kdf,\\
\L_{hK}(f_0df_1\wedge \dots \wedge df_p)
     &= \L_{hK}f_0.df_1\wedge \dots \wedge df_p)\\
&\quad + \sum_i(-1)^{i-1}f_0df_1\wedge \dots \L_{hK}df_i
     \dots \wedge df_p \otimes hf_i \\
&= \L_K(f_0df_1\wedge \dots \wedge df_p).
\endalign$$

\therosteritem7 can be seen as follows, using \therosteritem3,
\therosteritem2, and \therosteritem5:
$$\align
[hL,hf] &= [h(\ps\otimes Y),hf] = [\ps\wedge hY-d\pi^*Y\wedge \rh\ps,hf] \\
&= \ps\wedge [hY,hf] -\L_{hf}\ps\wedge hY - 0 \\
&\quad - d\pi^*Y\wedge [\rh\ps,hf] + \L_{hf}d\pi^*Y\wedge \rh\ps +0 \\
&= \ps\wedge h\L_Yf - di_Ydf\wedge \rh\ps
= h(\ps\wedge \L_Yf) = h\L_Lf
\endalign$$

\therosteritem8 We start with the following computation, using
\therosteritem4, \therosteritem5, and $i_{hX}\ps = i_X\ps $.
$$\align
\L_{hK}\pi^*L &= i_{hK}d\pi^*L - (-1)^{k-1}di_{hK}\pi^*L \\
&= i\left(-d\pi^*X\wedge\rh\ph + \ph \wedge hX\right)((-1)^l\ps\wedge
d\pi^*Y)\\
&\quad +(-1)^k d\;i\left(-d\pi^*X\wedge\rh\ph +
     \ph \wedge hX\right)(\ps\wedge \pi^*Y)\\
&=  -(-1)^l d\pi^*X\wedge i_{\rh\ph}\ps\wedge d\pi^*Y
     - (-1)^{l+(k-2)l}d\pi^*X\wedge \ps\wedge i_{\rh\ph}d\pi^*Y \\
&\quad +(-1)^l \ph \wedge i_{hX}\ps\wedge d\pi^*Y
     + \ph \wedge \ps\wedge i_{hX}d\pi^*Y\\
&\quad +(-1)^k d\left(-d\pi^*X\wedge i_{\rh\ph}\ps\wedge \pi^*Y
     +\ph \wedge i_{hX}\ps\wedge \pi^*Y\right)\\
&=  0 + (-1)^{(k-1)l}d\pi^*X\wedge \ps\wedge i_Y\ph
     +(-1)^l \ph \wedge i_{X}\ps\wedge d\pi^*Y \\
&\quad + \ph \wedge \ps\wedge \pi^*[X,Y]
     +0+ \ph \wedge di_{X}\ps\wedge \pi^*Y
     -(-1)^l \ph \wedge i_{X}\ps\wedge d\pi^*Y\\
&= \ph \wedge \ps\wedge \pi^*[X,Y]
     + \ph \wedge di_{X}\ps\wedge \pi^*Y
     -(-1)^{k+l} i_Y\ph \wedge \ps \wedge d\pi^*X,\\
\endalign$$
where we also used
$$\align
\L_{hX}\pi^*Y &= i_{hX}d\pi^*Y \\
&= i\left(-\sum \tfrac{\partial X^i}{\partial q^m}
     p_i\tfrac{\partial}{\partial p_m}
     + \sum X^i \tfrac{\partial}{\partial q^i}\right)
     \left(\sum \tfrac{\partial Y^k}{\partial q^n}p_kdq^n
     + \sum Y^kdp_k\right)\\
&= \sum\left( X^i\tfrac{\partial Y^k}{\partial q^i} -
     Y^i\tfrac{\partial X^k}{\partial q^i}\right)p_k = \pi^*[X,Y].
\endalign$$
Then we get
$$\align
\pi^*[K,L] &= \pi^*[\ph\otimes X,\ps\otimes Y],
     \quad\text{ use now \nmb!{2.1}.\thetag2} \\
&= \ph\wedge \ps\wedge\pi^*[X,Y] + \ph\wedge \L_X\ps\wedge \pi^*Y
     - \L_Y\ph\wedge \ps\wedge \pi^*X +0+0, \\
\L_{hK}\pi^*L &- \pi^*[K,L] = di_Y\ph\wedge \ps\wedge \pi^*X
     -(-1)^{k+l}i_Y\ph\wedge \ps\wedge d\pi^*X \\
&= d((-1)^{(k-1)l}\ps\wedge i_Y\ph\wedge \pi^*X)
     = (-1)^{(k-1)l}d\pi^*(i_LK).
\endalign$$

\therosteritem9 follows from \therosteritem8.

\therosteritem{10}
$i(\rh\pi^*K)\ps = i(\pi^*X.\rh\ph)\ps =
     \pi^*X.i_{\rh\ph}\ps=0$.

\therosteritem{11} We compute in turn
$$\align
\L_{\rh\pi^*K}\pi^*Y &= i(\pi^*X\wedge \rh\ph)d\pi^*Y =
     \pi^*X\wedge \L_{\rh\ph}\pi^*Y = - \pi^*X\wedge i_Y\ph =
     -\pi^*i_YK\\
\L_{\rh\pi^*K}\pi^*L &= \L_{\rh\pi^*K}(\pi^*Y\wedge \ps) =
     \L_{\rh\pi^*K}\pi^*Y\wedge \ps + \pi^*Y\wedge\L_{\rh\pi^*K}\ps\\
&= -\pi^*(i_YK)\wedge \ps + 0 \\
&= -(-1)^{(k-1)l}\pi^*i_LK
\endalign$$

\therosteritem{12} We have in turn
$$\align
i_{hX}\pi^*L &= i\left(-\sum \tfrac{\partial X^i}{\partial q^m}
     p_i\tfrac{\partial}{\partial p_m}
     + \sum X^i\tfrac{\partial}{\partial q^i}\right)(\ps\wedge \pi^*Y)\\
&= i_X\ps\wedge \pi^*Y = \pi^*i_XL,\\
i_{hK}\pi^*L &= i\left(-d\pi^*X\wedge \rh\ph
     + \ph\wedge hX\right)(\ps\wedge \pi^*Y)\\
&= -d\pi^*X\wedge i_{\rh\ph}\ps\wedge \pi^*Y
     + \ph\wedge i_{hX}\ps\wedge \pi^*Y\\
&= 0  + \ph\wedge i_{X}\ps\wedge \pi^*Y = \pi^*i_KL.
\endalign$$

\therosteritem{13}
 From \therosteritem8 we get
$$\align
i_{hK}d\pi^*L &= \L_{hK}\pi^*L +(-1)^{k-1}di_{hK}\pi^*L\\
&= \pi^*[K,L] -(-1)^kd\pi^*(i_KL+(-1)^{(k-1)(l-1)}i_LK).
\endalign$$

\therosteritem{14}
We just compute
$$\align
i_{hX}\rh(f_0df_1\wedge \dots\wedge df_q)
&= i_{hX}\left(\sum(-1)^{j-1}f_0df_1\wedge
     \dots\wedge hf_j\wedge \dots\wedge df_q \right)\\
&= \sum_{k<j}(-1)^{k+j}f_0df_1\wedge\dots i_{hX}df_k\wedge
     \dots\wedge hf_j\wedge \dots\wedge df_q \\
&\quad+ \sum_{k>j}(-1)^{k+j-1}f_0df_1\wedge
     \dots\wedge hf_j\wedge \dotsi_{hX}df_k\wedge\dots\wedge df_q \\
&=-\rh i_X(f_0df_1\wedge \dots\wedge df_q).
\endalign$$

\therosteritem{15} This is an easy consequence of \therosteritem4 and
\therosteritem5, namely
$$\align
\L_{h\ph}\pi^*L &= \L_{h\ph}(\ps\wedge\pi^*Y)
     = 0 + (-1)^{pl}\ps\wedge \L_{h\ph}^pi^*Y \\
&= - (-1)^{pl}\ps\wedge i_Y\ph = -(-1)^{pl}i_Ld\ph.
\endalign$$

\therosteritem{16} This can be seen by summing the following
evaluations.
$$\align
[\rh\pi^*K,h\ps] &= [\rh\pi^*(\ph\otimes X),h\ps] =
     [\rh\ph.\pi^*X,h\ps]\\
&= \pi^*X\wedge [\rh\ph,h\ps]
     -(-1)^{(k-1)q}\L_{h\ps}\pi^*X\wedge \rh\ph +
     (-1)^{k-1}d\pi^*X\wedge i_{\rh\ph}h\ps\\
&= \rh\ph\wedge i_Xd\ps,\\
\rh(i_Kd\ps) &= \rh(\ph\wedge i_Xd\ps) =
     \rh\ph\wedge i_Xd\ps +(-1)^k\ph\wedge \rh i_Xd\ps,\\
i_{hK}h\ps &= i(-d\pi^*X\wedge \rh\ph + \ph\wedge hX)(\rh d\ps)\\
&= -d\pi^*X\wedge i_{\rh\ph}\rh d\ps + \ph\wedge i_{hX}\rh d\ps = 0 -
     \ph\wedge i_{hX}\rh d\ps.\qed
\endalign$$
\enddemo

\subhead\nmb.{3.4} The extended insertion\endsubhead
For $A\in\Om^k(M;S^lTM)$ we define now the insertion operator
$$\gather
i_A:\Om^p(M;S^mTM)\to \Om^{p+k-1}(M;S^{m+l-1}TM)\\
i(\ph\otimes X_1\vee\dots\vee X_k)(\ps\otimes V) =
\ph\wedge \sum_ji_{X_j}\ps\otimes
     X_1\vee\dots\widehat{X_j}\dots\vee X_k\vee V.
\endgather$$
This is a graded derivation of degree $k-1$ of the graded commutative
algebra $\bigoplus_{m\ge0}\Om^m(M,STM)$ which vanishes on the
subalgebra $\Ga(STM)$.

\proclaim{Lemma. More formulas} For $A\in\Om^k(M;S^lTM)$, where
$l>0$, and $\ps\in \Om^q(M)$ we have on $T^*M$
\roster
\item $\L_{h\ps}\pi^*A = -(-1)^{qk}\pi^*i_Ad\ps$.
\item $[\rh\pi^*A,h\ps]=\rh\pi^*i_Ad\ps-(-1)^ki_{hA}h\ps$.
\endroster
\endproclaim

\demo{Proof}
\therosteritem1
We prove this by induction on $l$. For $l=1$ this is
\nmb!{3.3}.\therosteritem{15}. For the induction we compute as follows:
$$\align
\L_{h\ps}\pi^*(X\wedge A) &= \L_{h\ps}\pi^*X\wedge \pi^*A
     + \pi^*X\wedge \L_{h\ps}\pi^*A \\
&= - i_Xd\ps\wedge \pi^*A -(-1)^{qk}\pi^*X\wedge \pi^*i_Ad\ps \\
&= -(-1)^{qk}\pi^*(A\wedge i_Xd\ps +X\wedge i_Ad\ps)
     = -(-1)^{qk}\pi^*i_{X\wedge A}d\ps.
\endalign$$

\therosteritem2
We use again induction on $l$. For $l=1$ this is
\nmb!{3.3}.\therosteritem{16}. The left hand side equals:
$$\align
[\rh\pi^*&(X\wedge A),h\ps] = [\pi^*X\wedge \rh\pi^*A,h\ps]\\
&= \pi^*X\wedge [\rh\pi^*A,h\ps]
     -(-1)^{(k-1)q}\L_{h\ps}\pi^*X\wedge \rh\pi^*A +
     (-1)^kd\pi^*X\wedge i_{\rh\pi^*A}h\ps\\
&= \pi^*X\wedge \rh\pi^*i_Ad\ps
     -(-1)^k\pi^*X\wedge \pi^*i_{hA}h\ps \\
&\quad +(-1)^{(k-1)q}i_Xd\ps\wedge \rh\pi^*A +
     (-1)^kd\pi^*X\wedge i_{\rh\pi^*A}h\ps\\
\endalign$$
For the right hand side we get:
$$\align
\rh\pi^*&i_{X\wedge A}d\ps-(-1)^ki_{h(X\wedge A)}h\ps =
     \rh\pi^*(X\wedge i_Ad\ps + A\wedge i_Xd\ps) \\
&\quad -(-1)^ki(hX\wedge \pi^*A -d\pi^*X\wedge \rh\pi^*A
     + \pi^*X\wedge hA)h\ps \\
&= \pi^*X\wedge \rh\pi^*i_Ad\ps +
     \rh\pi^*A\wedge i_Xd\ps -(-1)^k\pi^*A\wedge\rh i_Xd\ps \\
&\quad -(-1)^k\pi^*A\wedge i_{hX}h\ps
     +(-1)^k d\pi^*X\wedge i_{\rh\pi^*A}h\ps
     -(-1)^k + \pi^*X\wedge i_{hA}h\ps.
\endalign$$
Using \nmb!{3.3}.\therosteritem{14} we see that it equals the left
hand side.
\qed\enddemo

\proclaim{\nmb.{3.5}. Theorem}
\roster
\item The linear injective mapping
$$ h:\Ga(\La T^*M\otimes TM) = \Om(M;TM) @>\pi^*>> \Om(T^*M) @>H>>
     \Om(T^*M;T(T^*M)) $$
     is a homomorphism for the Fr\"olicher-Nijenhuis brackets.
\item The linear mapping
$$ h:\Ga(STM) @>\pi^*>> \Om(T^*M) @>H>>
     \Om(T^*M;T(T^*M)) $$
     is a homomorphism from the symmetric Schouten bracket to the
     Fr\"olicher-Nijenhuis bracket. The kernel of $h$ is $H^0(M)$.
\item For differential forms $\ph,\ps\in \Om(M)$ we have
$$[h\ph,h\ps]=0.$$
\item For $A\in \Om(M;STM)$ and $\ps\in\Om(M)$ we have
$$\align
[hA,h\ps]&=hi_Ad\ps,\quad\text{ where}\\
i(\ph\otimes X_1\vee\dots\vee X_k)\ps:&=
     \ph\wedge \left(\sum_ji_{X_j}\ps\otimes
     X_1\vee\dots\widehat{X_j}\dots\vee X_k\right).
\endalign$$
\item For $\dim M\ge2$, in general
       $[h\Om^{k_1}(M;S^{l_1}TM),h\Om^{k_2}(M;S^{l_2}TM)]$
       does not lie in the image of $h$, if $k_1,l_1\ge1$ and
       $l_2\ge2$ (or under the symmetric condition).
\endroster
\endproclaim

\demo{Proof} \therosteritem1
We have to show that $[hK,hL]=h[K,L]$ for $K\in\Om^k(M;TM)$ and
$L\in\Om^l(M;TM)$ and we do this by induction on $k+l$. The case of
vector fields $k+l=0$ is well know, see \nmb!{3.3}.\therosteritem1.
Since the question is local on $M$ and since $\Om^{k+1}(M;TM)$ is
locally linearly generated by $df\wedge K$ for $f\in\Om^0(M)$ and
$K\in\Om^k(M;TM)$ it suffices to check that $[hK,hL]=h[K,L]$ implies
$[h(df\wedge K),hL]=h[df\wedge K,L]$.
We have
$$h(df\wedge K)= \rh d(df\wedge \pi^*K) = -d\pi^*K\wedge hf +
df\wedge hK.$$
Using twice \nmb!{2.1}.\thetag4 we get then
$$\align
[h(df\wedge K),hL] &=
     df\wedge [hK,hL] -(-1)^{(1+k)l}\L_{hL}df\wedge hK +0 \\
&\quad -d\pi^*K\wedge [hf,hL] +(-1)^{(1+k)l}\L_{hL}d\pi^*K\wedge hf-0\\
&= df\wedge h[K,L] -(-1)^{(1+k)l}\L_{L}df\wedge hK \\
&\quad +d\pi^*K\wedge h\L_Lf +(-1)^{kl}d\pi^*[L,K]\wedge hf,
\endalign$$
where we used in turn induction, \nmb!{3.3}.\thetag6,
\nmb!{3.3}.\thetag7, and \nmb!{3.3}.\thetag8.
On the other hand we have again by \nmb!{2.1}.\thetag4
$$\align
h[df\wedge K,L] &= \rh d\pi^*\left(df\wedge [K,L]
     -(-1)^{(1+k)l}\L_Ldf\wedge K +0 \right)\\
&= -hf\wedge d\pi^*[K,L] + df\wedge h[K,L] \\
&\quad +(-1)^{kl}\rh\L_Ldf\wedge d\pi^*K +(-1)^{kl+l+1}\L_Ldf\wedge hK,
\endalign$$
which equals the expression for $[h(df\wedge K),hL]$ from above.

\therosteritem2 This is well known, and easy to check starting from
\nmb!{3.3}.\thetag1.

\therosteritem3 is \nmb!{3.3}.\therosteritem2.

\therosteritem4
First we prove a partial result.\newline
{\bf Claim.} For $K\in \Om^k(M;TM)$ and $\ps\in\Om^q(M)$ we have
$[hK,h\ps]=h\L_K\ps = hi_K\ps$.

To check the claim we use induction on $q=\operatorname{deg}\ps$.
For $q=0$ this is \nmb!{3.3}.\thetag7. Since the assertion is local on
$M$ it suffices to consider $df\wedge \ps$ for the induction step. In
the following computation we use \nmb!{2.1}.\thetag4, induction,
\nmb!{3.3}.\therosteritem6, and \nmb!{3.3}.\therosteritem7:
$$\align
[hK,h&(df\wedge \ps)] = [hK, df\wedge h\ps - d\ps \wedge hf]\\
&= (-1)^kdf\wedge [hK,h\ps] + \L_{hK}df\wedge h\ps - 0\\
&\quad -(-1)^{(q+1)k} d\ps\wedge [hK,hf] - \L_{hK}d\ps\wedge hf +0\\
&= (-1)^kdf\wedge h\L_K\ps + h\L_Kdf\wedge h\ps
     -(-1)^{(q+1)k} d\ps\wedge h\L_Kf - \L_Kd\ps\wedge hf \\
\endalign$$
On the other hand we have by \nmb!{2.1}.\thetag3
$$\align
h\L_K&(df\wedge \ps) = h\left((-1)^kd\L_Kf\wedge \ps +
     (-1)^kdf\wedge \L_K\ps  \right)\\
&= \rh\left(d\L_Kf\wedge d\ps -
     (-1)^kdf\wedge d\L_K\ps  \right)\\
&= -h\L_Kf\wedge d\ps +(-1)^k
d\L_Kf\wedge h\ps -(-1)^khf\wedge d\L_K\ps +(-1)^kdf\wedge hL_K\ps,
\endalign$$
which equals the above expression. So the claim follows.

Now we can extend this result to $A\in\Om^k(M;S^lTM)$ by induction on
$l$. For $l=1$ this is the claim above.
For the induction we compute first the left hand side, using
\nmb!{2.1}.\thetag4, the claim, \nmb!{3.4}, induction,
and \nmb!{3.3}.\therosteritem{14}
$$\align
[h(X.A),h\ps] &= [hX\wedge \pi^*A -d\pi^*X\wedge \rh\pi^*A +
     \pi^*X\wedge hA,h\ps]\\
&= \pi^*A\wedge [hX,h\ps] -(-1)^{kq} \L_{h\ps}\pi^*A\wedge hX
     +(-1)^k d\pi^*A\wedge i_{hX}h\ps\\
&\quad- d\pi^*X\wedge [\rh\pi^*A,h\ps]
     +(-1)^{kq}\L_{h\ps}d\pi^*X\wedge\rh\pi^*A - 0 \\
&\quad+ \pi^*X\wedge [hA,h\ps] -(-1)^{kq} \L_{h\ps}\pi^*X\wedge hA
     +(-1)^k d\pi^*X\wedge i_{hA}h\ps\\
&= \pi^*A\wedge hi_Xd\ps + \pi^*i_Ad\ps\wedge hX
     -(-1)^k d\pi^*A\wedge \rh i_Xd\ps\\
&\quad - d\pi^*X\wedge\rh i_Ad\ps
     -(-1)^{(k-1)q}di_Xd\ps\wedge\rh\pi^*A \\
&\quad + \pi^*X\wedge hi_Ad\ps +(-1)^{kq} i_Xd\ps\wedge hA.
\endalign$$
The right hand side is
$$\align
hi(&X\wedge A)d\ps = h(A\wedge i_Xd\ps + X\wedge i_Ad\ps)\\
&= \rh(d\pi^*A\wedge i_Xd\ps +(-1)^k\pi^*A\wedge di_Xd\ps
     + d\pi^*X\wedge \pi^*i_Ad\ps + \pi^*X\wedge d\pi^*i_Ad\ps)\\
&= hA\wedge i_Xd\ps +(-1)^kd\pi^*A\wedge \rh i_Xd\ps
     +(-1)^k\rh\pi^*A\wedge di_Xd\ps + \pi^*A\wedge hi_Xd\ps \\
&\quad + hX\wedge \pi^*i_Ad\ps - d\pi^*X\wedge \rh\pi^*i_Ad\ps
     + \pi^*X\wedge h\pi^*i_Ad\ps),
\endalign$$
which equals the left hand side.

\therosteritem5
By \nmb!{3.2}.\therosteritem3 the image of $\pi^*:\Om(M;S^lTM)$ is
the space of all horizontal forms $\Ph\in\Om(T^*M)$ satisfying
$\L_I\Ph=l.\Ph$. In local coordinates on $M$ we consider then, using
the bracket $\{\quad,\quad\}^1$ described in \nmb!{2.3},
$$\align
\{\pi^*(dq^1\otimes \frac{\partial}{\partial q^1}),
     \pi^*(\frac{\partial}{\partial q^1}\frac{\partial}{\partial q^2})\}^1 &=
     \{p_1dq^1,p_1p_2\}^1 = \{p_1,p_1p_2\}dq^1 - \{q^1,p_1p_2\}dp_1\\
&= p_2dp_1,\\
d\{p_1dq^1,p_1p_2\}^1 &= -dp_1\wedge dp_2.
\endalign$$
Thus $\{p_1dq^1,p_1p_2\}^1$ plus something exact can never be
horizontal.
\qed\enddemo

\enddocument

\bye